\documentclass{mrlart7}
  \def\volno{18}
  \def\yrno{2011}
  \def\issueno{00}
  \def\lpageno{100NN} 
\usepackage{amsmath}
\usepackage{amsxtra}
\usepackage{amscd}
\usepackage{amsthm}
\usepackage{amsfonts}
\usepackage{amssymb}
\usepackage{eucal}
\usepackage[matrix,arrow,curve]{xy}

\theoremstyle{plain}
\newtheorem{Thm}[subsection]{Theorem}
\newtheorem{Cor}[subsection]{Corollary}
\newtheorem{Lem}[subsection]{Lemma}
\newtheorem{Prop}[subsection]{Proposition}
\newtheorem{Conj}[subsection]{Conjecture}

\theoremstyle{definition}
\newtheorem{Def}[subsection]{Definition}

\theoremstyle{remark}

\newtheorem{Rem}[subsection]{Remark}

\newtheorem{conj}{Conjecture}

\numberwithin{equation}{section}
\renewcommand{\rm}{\normalshape}

\newif\ifShowLabels
\ShowLabelstrue
\newdimen\theight
\def\TeXref#1{%
    \leavevmode\vadjust{\setbox0=\hbox{{\tt
        \quad\quad  {\small \rm #1}}}%
    \theight=\ht0
    \advance\theight by \lineskip
    \kern -\theight \vbox to
    \theight{\rightline{\rlap{\box0}}%
    \vss}%
    }}%

\ShowLabelsfalse

\renewcommand{\sec}[2]{\section{#2}\label{S:#1}%
    \ifShowLabels \TeXref{{S:#1}} \fi}
\newcommand{\ssec}[2]{\subsection{#2}\label{SS:#1}%
    \ifShowLabels \TeXref{{SS:#1}} \fi}

\newcommand{\refs}[1]{Section ~\ref{S:#1}}
\newcommand{\refss}[1]{Section ~\ref{SS:#1}}

\newcommand{\reft}[1]{Theorem ~\ref{T:#1}}
\newcommand{\refl}[1]{Lemma ~\ref{L:#1}}

\newcommand{\refe}[1]{\eqref{E:#1}}
\newcommand{\refco}[1]{Conjecture ~\ref{Co:#1}}

\newenvironment{thm}[1]%
    { \begin{Thm} \label{T:#1}  \ifShowLabels \TeXref{T:#1} \fi }%
    { \end{Thm} }

\renewcommand{\th}[1]{\begin{thm}{#1} \sl }
\renewcommand{\eth}{\end{thm} }

\newenvironment{lemma}[1]%
    { \begin{Lem} \label{L:#1}  \ifShowLabels \TeXref{L:#1} \fi }%
    { \end{Lem} }
\newcommand{\lem}[1]{\begin{lemma}{#1} \sl}
\newcommand{\elem}{\end{lemma}}

\newenvironment{propos}[1]%
    { \begin{Prop} \label{P:#1}  \ifShowLabels \TeXref{P:#1} \fi }%
    { \end{Prop} }
\newcommand{\prop}[1]{\begin{propos}{#1}\sl }
\newcommand{\eprop}{\end{propos}}

\newenvironment{corol}[1]%
    { \begin{Cor} \label{C:#1}  \ifShowLabels \TeXref{C:#1} \fi }%
    { \end{Cor} }
\newcommand{\cor}[1]{\begin{corol}{#1} \sl }
\newcommand{\ecor}{\end{corol}}

\newenvironment{defeni}[1]%
    { \begin{Def} \label{D:#1}  \ifShowLabels \TeXref{D:#1} \fi }%
    { \end{Def} }
\newcommand{\defe}[1]{\begin{defeni}{#1} \sl }
\newcommand{\edefe}{\end{defeni}}

\newenvironment{remark}[1]%
    { \begin{Rem} \label{R:#1}  \ifShowLabels \TeXref{R:#1} \fi }%
    { \end{Rem} }
\newcommand{\rem}[1]{\begin{remark}{#1}}
\newcommand{\erem}{\end{remark}}

\newenvironment{conjec}[1]%
    { \begin{Conj} \label{Co:#1}  \ifShowLabels \TeXref{Co:#1} \fi }%
    { \end{Conj} }
\renewcommand{\conj}[1]{\begin{conjec}{#1} \sl }
\newcommand{\econj}{\end{conjec}}

\newcommand{\eq}[1]%
    { \ifShowLabels \TeXref{E:#1} \fi
       \begin{equation} \label{E:#1} }
\newcommand{\eeq}{ \end{equation} }

\newcommand{\prf}{ \begin{proof} }
\newcommand{\epr}{ \end{proof} }

\newcommand\alp{\alpha}

\newcommand\lam{\lambda}        \newcommand\Lam{\Lambda}

\newcommand\calW{{\mathcal{W}}}

\newcommand\GG{\mathbb{G}}

\newcommand\ZZ{\mathbb{Z}}

 \newcommand\grg{{\mathfrak{g}}}

\newcommand\sdp{\times \hskip -0.3em {\raise 0.3ex
\hbox{$\scriptscriptstyle |$}}} 

\newcommand\Gr{\operatorname{Gr}}

\newcommand\Hom{\operatorname {Hom}}

\newcommand\Int{\operatorname{Int}}

\newcommand\olam{{\overline{\lambda}}}

\newcommand\hatG{{\widehat{G}}}

\newcommand\tilG{{\widetilde{G}}}

\newcommand\tilS{{\widetilde{S}}}

\newcommand\x{\times}

\newcommand{\la}{\langle}

\newcommand\nc{\newcommand}

\newcommand{\IC}{{\operatorname{IC}}}

\newcommand{\iso}{{\stackrel{\sim}{\longrightarrow}}}

\nc\aff{\operatorname{aff}}
\nc\oGr{\overline{\Gr}}
\nc\Bun{\operatorname{Bun}}
\nc\hgrg{\widehat{\grg}}
\renewcommand\Int{\operatorname{Int}}
\nc\bInt{\overline{\Int}}
\nc\hatLam{\widehat{\Lam}}
\nc\bmu{\overline{\mu}}
\nc\bnu{\overline{\nu}}
\nc\blambda{\overline{\lam}}

\nc\ocalW{\overline{\calW}}
\nc\pos{\operatorname{pos}}
\nc\IH{\operatorname{IH}}
\nc\Rep{\operatorname{Rep}}
\nc\Gal{\operatorname{Gal}}
\nc{\tilGr}{\widetilde{\Gr}}

\nc\Pic{\operatorname{Pic}}

\nc{\HC}{{\mathcal{HC}}}
\nc{\on}{\operatorname}
\nc{\BA}{{\mathbb{A}}}
\nc{\BC}{{\mathbb{C}}}
\nc{\BM}{{\mathbb{M}}}
\nc{\BN}{{\mathbb{N}}}
\nc{\BP}{{\mathbb{P}}}
\nc{\BR}{{\mathbb{R}}}
\nc{\BZ}{{\mathbb{Z}}}
\nc{\BS}{{\mathbb{S}}}

\nc{\CA}{{\mathcal{A}}}
\nc{\CB}{{\mathcal{B}}}
\nc{\CalD}{{\mathcal D}}
\nc{\CE}{{\mathcal{E}}}
\nc{\CF}{{\mathcal{F}}}
\nc{\CG}{{\mathcal{G}}}
\nc{\CH}{{\mathcal{H}}}
\nc{\CK}{{\mathcal{K}}}
\nc{\CL}{{\mathcal{L}}}
\nc{\CM}{{\mathcal{M}}}
\nc{\CMM}{{\mathcal{M}^{\operatorname{gen}}_\hbar(-\rho)}}
\nc{\CN}{{\mathcal{N}}}
\nc{\CO}{{\mathcal{O}}}
\nc{\CP}{{\mathcal{P}}}
\nc{\CQ}{{\mathcal{Q}}}
\nc{\CR}{{\mathcal{R}}}
\nc{\CS}{{\mathcal{S}}}
\nc{\CT}{{\mathcal{T}}}
\nc{\CU}{{\mathcal{U}}}
\nc{\CV}{{\mathcal{V}}}
\nc{\CW}{{\mathcal{W}}}
\nc{\CX}{{\mathcal{X}}}
\nc{\CZ}{{\mathcal{Z}}}

\nc{\gen}{{\operatorname{gen}}}
\nc{\cM}{{\check{\mathcal M}}{}}
\nc{\csM}{{\check{\mathcal A}}{}}
\nc{\obM}{{\overset{\circ}{\mathbf M}}{}}
\nc{\oCA}{{\overset{\circ}{\mathcal A}}{}}
\nc{\obA}{{\overset{\circ}{\mathbf A}}{}}
\nc{\ooM}{{\overset{\circ}{M}}{}}
\nc{\osM}{{\overset{\circ}{\mathsf M}}{}}
\nc{\vM}{{\overset{\bullet}{\mathcal M}}{}}
\nc{\nM}{{\underset{\bullet}{\mathcal M}}{}}
\nc{\obD}{{\overset{\circ}{\mathbf D}}{}}
\nc{\cp}{{\overset{\circ}{\mathbf p}}{}}
\nc{\ofZ}{{\overset{\circ}{\mathfrak Z}}{}}

\nc{\fa}{{\mathfrak{a}}}
\nc{\fb}{{\mathfrak{b}}}
\nc{\fg}{{\mathfrak{g}}}
\nc{\fgl}{{\mathfrak{gl}}}
\nc{\fh}{{\mathfrak{h}}}
\nc{\fri}{{\mathfrak{i}}}
\nc{\fj}{{\mathfrak{j}}}
\nc{\fm}{{\mathfrak{m}}}
\nc{\fn}{{\mathfrak{n}}}
\nc{\ft}{{\mathfrak{t}}}
\nc{\fu}{{\mathfrak{u}}}
\nc{\fp}{{\mathfrak{p}}}
\nc{\frr}{{\mathfrak{r}}}
\nc{\fs}{{\mathfrak{s}}}
\nc{\fT}{{\mathfrak{T}}}
\nc{\ofT}{{\overline{\mathfrak T}}}
\nc{\ofS}{{\overline{\mathfrak S}}}
\nc{\fsl}{{\mathfrak{sl}}}
\nc{\hsl}{{\widehat{\mathfrak{sl}}}}
\nc{\hgl}{{\widehat{\mathfrak{gl}}}}
\nc{\hg}{{\widehat{\mathfrak{g}}}}
\nc{\chg}{{\widehat{\mathfrak{g}}}{}^\vee}
\nc{\hn}{{\widehat{\mathfrak{n}}}}
\nc{\chn}{{\widehat{\mathfrak{n}}}{}^\vee}

\nc{\fA}{{\mathfrak{A}}}
\nc{\fB}{{\mathfrak{B}}}
\nc{\fD}{{\mathfrak{D}}}
\nc{\fE}{{\mathfrak{E}}}
\nc{\fF}{{\mathfrak{F}}}
\nc{\fG}{{\mathfrak{G}}}
\nc{\fK}{{\mathfrak{K}}}
\nc{\fL}{{\mathfrak{L}}}
\nc{\fM}{{\mathfrak{M}}}
\nc{\fN}{{\mathfrak{N}}}
\nc{\frP}{{\mathfrak{P}}}
\nc{\fS}{{\mathfrak S}}
\nc{\fU}{{\mathfrak{U}}}
\nc{\fZ}{{\mathfrak{Z}}}

\nc{\bb}{{\mathbf{b}}}
\nc{\bc}{{\mathbf{c}}}
\nc{\be}{{\mathbf{e}}}
\nc{\bj}{{\mathbf{j}}}
\nc{\bn}{{\mathbf{n}}}
\nc{\bp}{{\mathbf{p}}}
\nc{\bq}{{\mathbf{q}}}
\nc{\bv}{{\mathbf{v}}}
\nc{\bx}{{\mathbf{x}}}
\nc{\by}{{\mathbf{y}}}
\nc{\bw}{{\mathbf{w}}}
\nc{\bA}{{\mathbf{A}}}
\nc{\bB}{{\mathbf{B}}}
\nc{\bC}{{\mathbf{C}}}
\nc{\bK}{{\mathbf{K}}}
\nc{\bD}{{\mathbf{D}}}
\nc{\bH}{{\mathbf{H}}}
\nc{\bM}{{\mathbf{M}}}
\nc{\bN}{{\mathbf{N}}}
\nc{\bS}{{\mathbf{S}}}
\nc{\bT}{{\mathbf{T}}}
\nc{\bV}{{\mathbf{V}}}
\nc{\bW}{{\mathbf{W}}}
\nc{\bX}{{\mathbf{X}}}
\nc{\bP}{{\mathbf{P}}}
\nc{\bZ}{{\mathbf{Z}}}

\nc{\sA}{{\mathsf{A}}}
\nc{\sB}{{\mathsf{B}}}
\nc{\sC}{{\mathsf{C}}}
\nc{\sD}{{\mathsf{D}}}
\nc{\sF}{{\mathsf{F}}}
\nc{\sK}{{\mathsf{K}}}
\nc{\sM}{{\mathsf{M}}}
\nc{\sO}{{\mathsf{O}}}
\nc{\sQ}{{\mathsf{Q}}}
\nc{\sP}{{\mathsf{P}}}
\nc{\sR}{{\mathsf{R}}}
\nc{\sV}{{\mathsf{V}}}
\nc{\sZ}{{\mathsf{Z}}}
\nc{\sfp}{{\mathsf{p}}}
\nc{\sr}{{\mathsf{r}}}
\nc{\sfb}{{\mathsf{b}}}
\nc{\sfc}{{\mathsf{c}}}
\nc{\sd}{{\mathsf{d}}}
\nc{\sfl}{{\mathsf{l}}}

\nc{\BK}{{\bar{K}}}

\nc{\tA}{{\widetilde{\mathbf{A}}}}
\nc{\tB}{{\widetilde{\mathcal{B}}}}
\nc{\tg}{{\widetilde{\mathfrak{g}}}}
\nc{\tG}{{\widetilde{G}}}
\nc{\TM}{{\widetilde{\mathbb{M}}}{}}
\nc{\tO}{{\widetilde{\mathsf{O}}}{}}
\nc{\tU}{{\widetilde{\mathfrak{U}}}{}}
\nc{\TZ}{{\tilde{Z}}}
\nc{\tx}{{\tilde{x}}}
\nc{\tbv}{{\tilde{\bv}}}
\nc{\tfP}{{\widetilde{\mathfrak{P}}}{}}
\nc{\tz}{{\tilde{\zeta}}}
\nc{\tmu}{{\tilde{\mu}}}

\nc{\td}{\ddot{\underline{d}}{}}
\nc{\tzeta}{\widetilde{\zeta}{}}
\nc{\hd}{{\widehat{\underline{d}}}}
\nc{\hG}{{\widehat{G}}}
\nc{\hBP}{\widehat{\mathbb P}{}}
\nc{\hQ}{{\widehat{Q}}}
\nc{\hsM}{\widehat{\mathsf M}{}}
\nc{\hfM}{\widehat{\mathfrak M}{}}
\nc{\hCP}{\widehat{\mathcal P}{}}
\nc{\hCR}{\widehat{\mathcal R}{}}
\nc{\hCS}{{\widehat{\mathcal S}}}
\nc{\hfZ}{\widehat{\mathfrak Z}{}}

\nc{\urho}{\underline{\rho}}
\nc{\uB}{\underline{B}}
\nc{\uC}{{\underline{\mathbb{C}}}}
\nc{\ui}{\underline{i}}
\nc{\ofP}{{\overline{\mathfrak{P}}}}

\nc{\hrho}{{\hat{\rho}}}

\nc{\unl}{\underline}
\nc{\ol}{\overline}
\nc{\one}{{\mathbf{1}}}
\nc{\two}{{\mathbf{t}}}

\nc{\Tot}{{\mathop{\operatorname{\rm Tot}}}}
\nc{\Hilb}{{\mathop{\operatorname{\rm Hilb}}}}
\nc{\CHom}{{\mathop{\operatorname{{\mathcal{H}}\it om}}}}
\nc{\defi}{{\mathop{\operatorname{\rm def}}}}
\nc{\length}{{\mathop{\operatorname{\rm length}}}}

\nc{\Cliff}{{\mathsf{Cliff}}}
\nc{\Fl}{{\mathsf{Fl}}}
\nc{\Fib}{{\mathsf{Fib}}}
\nc{\Coh}{{\mathsf{Coh}}}
\nc{\FCoh}{{\mathsf{FCoh}}}

\nc{\reg}{{\text{\rm reg}}}

\nc{\cplus}{{\mathbf{C}_+}}
\nc{\cminus}{{\mathbf{C}_-}}
\nc{\cthree}{{\mathbf{C}_*}}
\nc{\Qbar}{{\bar{Q}}}

\nc{\bh}{{\bar{h}}}
\nc{\bOmega}{{\overline{\Omega}}}
\nc\tGr{\widetilde{\Gr}}

\nc{\seq}[1]{\stackrel{#1}{\sim}}

\nc\uS{\underline{S}}
\begin{document}
\title[Dynamical Weyl groups and the affine Grassmannian]{Dynamical Weyl groups and equivariant cohomology of transversal slices on affine Grassmannians}
\author{Alexander Braverman and Michael Finkelberg}

\begin{abstract}Let $G$ be a reductive group and let $\check{G}$ be its Langlands dual.
We give an interpretation of the dynamical Weyl group of $\check{G}$ defined in \cite{EV} in terms of the
geometry of the affine Grassmannian $\Gr$ of $G$.
In this interpretation the dynamical parameters of \cite{EV} correspond to equivariant parameters with respect to
certain natural torus acting on $\Gr$. We also present a conjectural generalization of our results to the case of affine
Kac-Moody groups.
\end{abstract}
\maketitle

\sec{notstat}{Introduction and statements of the results}

\ssec{}{Notations and overview}

We mainly follow the notations of~\cite{BF}. So $G\supset B\supset T$ is a reductive complex algebraic group with
a Borel subgroup, and a Cartan subgroup. Let $\check G\supset\check T$ denote the
Langlands dual group with its Cartan torus. Let $\CK=\BC((t)),\ \CO=\BC[[t]]$. The affine Grassmannian
$\Gr=\Gr_G=G(\CK)/G(\CO)$ carries the category $\on{Perv}_{G(\CO)}(\Gr)$ of $G(\CO)$-equivariant perverse constructible sheaves. It is equipped with the convolution monoidal structure, and is tensor equivalent to the tensor category $\on{Rep}(\check G)$ of representations of $\check G$ (see~\cite{Lu-qan},~\cite{Gi},~\cite{MV},~\cite{BD}). We denote by $\tilde S:\ \on{Rep}(\check G)\to\on{Perv}_{G(\CO)}(\Gr)$ the geometric Satake isomorphism functor, and by $S:\ \on{Rep}(\check G)\to\on{Perv}_{G(\CO)\rtimes\BC^*}(\Gr)$ its extension to the monoidal category of $G(\CO)\rtimes\BC^*$-equivariant perverse constructible sheaves.

The Lie algebras of $\check G\supset\check T$ are denoted by $\check\fg\supset\check\ft$. We have
$\check\ft=\ft^*$ canonically. The Weyl group of
$G,\check G$ is denoted by $W$. Let $\Lambda=\Lambda_G$ denote the coweight lattice of $T$ which is the same as the weight lattice of $\check T$. The choice of Borel subgroup $B\subset G$ defines the cone $\Lambda^+_G\subset\Lambda_G$ of dominant weights of $\check T$.

The lattice $\Lambda=\Lambda_G$ is identified with the quotient $T(\CK)/T(\CO)$. For $\lambda\in\Lambda$
we denote by $t^\lambda$ any lift of $\lambda$ to $T(\CK)$. Its image in $\Gr_G=G(\CK)/G(\CO)$ is independent
of the choice of a lift, and will be also denoted by $t^\lambda$, or sometimes just $\lambda\in\Gr_G$.
Moreover, we will keep the same name for the closed embedding $\lambda\hookrightarrow\Gr$.
Let $\Gr^\lambda$ denote the $G(\CO)$-orbit of $\lambda$, and let $\ol{\Gr}{}^\lambda\subset\Gr$ denote the
closure of $\Gr^\lambda$. It is well known that $\Gr=\bigcup_{\lambda\in\Lambda}\Gr^\lambda$, and that
$\Gr^\lambda=\Gr^\mu$ iff $\lambda$ and $\mu$ lie in the same $W$-orbit on $\Lambda$. In particular,
$\Gr=\bigsqcup_{\lambda\in\Lambda^+}\Gr^\lambda$.

It is known that as a byproduct of the existence of $S$ (or $\tilS$) one can get a geometric construction of finite-dimensional
representations of $\check{G}$. It is not difficult to note that this realization immediately gives rise to some additional structures
on finite-dimensional representations of $\check{G}$ and one may wonder what is the representation-theoretic
meaning of those finer structures. In this paper we are going to explore one such example. Namely,
let $V\in \Rep(\check{G})$ and let $\mu$ be a weight of $V$, i.e.
an element of $\Lam$ such that $V_{\mu}\neq 0$.
Let also $w\in W$. Then we explain in \refss{hyp} that from the geometry of $\Gr$ one gets a canonical isomorphism
$V_{\mu}\otimes\BC(\ft\times\BA^1)\to
V_{w\mu}\otimes\BC(\ft\times\BA^1)$. Here $\BA^1$ is the affine line and $\BC(\ft\times\BA^1)$ denotes
the field of rational functions on $\ft\times\BA^1$. The main purpose of this paper is to show that
when $\mu\in\Lam^+$ the above isomorphism coincides (up to some simple reparametrization) with the one
provided by the {\em dynamical Weyl group} constructed in \cite{EV}.

\ssec{}{More notation}For $\lambda,\mu\in\Lambda^+$, we have $\ol{\Gr}{}^\mu\subset\ol{\Gr}{}^\lambda$ iff $\mu\leq\lambda$, i.e.
$\lambda-\mu$ is a sum of positive roots of $\check G$. We denote by $G[t^{-1}]_1\subset G[t^{-1}]\subset G(\CK)$ the kernel of the natural (``evaluation at $\infty$") homomorphism $G[t^{-1}]\to G$. For $\mu\in\Lambda$ we set $\Gr_\mu=G[t^{-1}]\cdot\mu\subset\Gr$, and $\CW_\mu=G[t^{-1}]_1\cdot\mu$. The locally closed embedding $\CW_\mu\hookrightarrow\Gr$ will be denoted by $\fri_\mu$. For dominant $\mu\leq\lambda\in\Lambda^+$, we set $\CW^\lambda_\mu=\Gr^\lambda\cap\CW_\mu$, and $\ol{\CW}{}^\lambda_\mu=
\ol{\Gr}{}^\lambda\cap\CW_\mu$. The variety $\ol{\CW}{}^\lambda_\mu$ can be thought of as a transversal slice
to $\Gr^\mu$ inside $\ol{\Gr}{}^\lambda$.

For $V\in\on{Rep}(\check G)$ we can restrict the $G(\CO)\rtimes\BC^*$-equivariant perverse sheaf $S(V)$
to the transversal slice $\CW_\mu$. We have $\fri_\mu^*S(V)[-\dim\Gr^\mu]\simeq\fri_\mu^!S(V)[\dim\Gr^\mu]$, and we will denote this restriction by $S_\mu(V)$. It is a $T\times\BC^*$-equivariant perverse sheaf on $\CW_\mu$.

Let $T\subset B_-\subset G$ be the opposite (to $T\subset B\subset G$) Borel subgroup, and let $N_-\subset B_-$ be its unipotent radical. For $w\in W$ we consider the conjugate subgroup $N_-^w:=\tilde{w}{}^{-1}N_-\tilde w$ for any lift $\tilde w$ of $w$ to the normalizer of $T$. For $\lambda\in\Lambda$ we set $\fT_\lambda^w=N_-^w(\CK)\cdot\lambda\subset\Gr$. The locally closed embedding $\fT_\lambda^w\hookrightarrow\Gr$ is denoted by $\imath_\lambda^w$, and the closed embedding $\lambda\hookrightarrow\fT_\lambda^w$ is denoted by
$\fj_\lambda^w$.

For a $T\times\BC^*$-equivariant constructible complex $\CF$ on $\Gr$, we denote by $R^w_\mu\CF$ the
$H^\bullet_{T\times\BC^*}(pt)$-module $\fj_\mu^{w*}\imath_\mu^{w!}\CF[\langle w\mu,2\check\rho\rangle]$.
Here $2\check\rho$ is the sum of positive roots of $G$, and $H^\bullet_{T\times\BC^*}(pt)\simeq\BC[\ft\times\BA^1]$. We denote by $r^w_\mu\CF$ the vector space
$\fj_\mu^{w*}\imath_\mu^{w!}\CF[\langle w\mu,2\check\rho\rangle]$ with $T\times\BC^*$-equivariance forgotten. Thus, $r^w_\mu\CF$ is the fiber of $R^w_\mu\CF$ over $0\in\ft\times\BA^1$.
Similarly, for a $T\times\BC^*$-equivariant constructible complex $\CF$ on $\CW_\mu$, we denote by $\sR^w_\mu\CF$ the $H^\bullet_{T\times\BC^*}(pt)$-module $\fj_\mu^{w*}\imath_\mu^{w!}\CF$.
We denote by $\sr^w_\mu\CF$ the vector space
$\fj_\mu^{w*}\imath_\mu^{w!}\CF$ with $T\times\BC^*$-equivariance forgotten. Thus, $\sr^w_\mu\CF$ is the fiber of $\sR^w_\mu\CF$ over $0\in\ft\times\BA^1$.

\ssec{hyp}{Hyperbolic stalks}

Recall that for $V\in\on{Rep}(\check G)$ we have a $T\times\BC^*$-equivariant perverse sheaf $S_\mu(V)$ on the slice $\CW_\mu$. The equivariant costalk $\mu^!S_\mu(V)=\fj_\mu^{w!}\imath_\mu^{w!}S_\mu(V)$ is an $H^\bullet_{T\times\BC^*}(pt)$-module equipped with a natural morphism $C_\mu^w$ to $\fj_\mu^{w*}\imath_\mu^{w!}S_\mu(V)=\sR^w_\mu S_\mu(V)$. Note that $C_\mu^w$ is a morphism of
$H^\bullet_{T\times\BC^*}(pt)=\BC[\ft\times\BA^1]$-modules. If we extend the scalars to the field of fractions $\BC(\ft\times\BA^1)$, the morphism $C_\mu^w$ becomes an isomorphism due to Localization Theorem in equivariant
cohomology (note that $\mu$ is the unique $T$-fixed point of $\CW_\mu$). For $w,y\in W$, we denote by
$\sA_{w,y}^{S_\mu(V)}:\ \sR^w_\mu S_\mu(V)\otimes_{\BC[\ft\times\BA^1]}\BC(\ft\times\BA^1)\iso
\sR^y_\mu S_\mu(V)\otimes_{\BC[\ft\times\BA^1]}\BC(\ft\times\BA^1)$ the composition
$C_\mu^y\circ(C_\mu^w)^{-1}$.

\lem{ff}
(a) For $w\in W,\ \mu\in\Lambda,\ V\in\on{Rep}(\check G)$, we have a canonical isomorphism
$\sr^w_\mu S_\mu(V)\cong V_{w\mu}$ (the $w\mu$-weight space of $V$);

(b) For $w\in W,\ \mu\in\Lambda,\ V\in\on{Rep}(\check G)$, we have a canonical isomorphism
$\sR^w_\mu S_\mu(V)\cong \sr^w_\mu S_\mu(V)\otimes\BC[\ft\times\BA^1]$.
\elem

According to~\refl{ff}, we can view the isomorphism $\sA_{w,y}^{S_\mu(V)}$ as going from
$V_{w\mu}\otimes\BC(\ft\times\BA^1)$ to $V_{y\mu}\otimes\BC(\ft\times\BA^1)$. Let us denote by $\hbar\in\BC[\BA^1]$ the positive generator of $H^2_{\BC^*}(pt,\BZ)$.

\ssec{dyn}{Dynamical Weyl group}
The dynamical Weyl group operators for simple Lie algebras were introduced in~\cite{Z} by analogy with extremal
projectors of~\cite{AST}. We will follow the presentation of~\cite{EV}.
We refer the reader to {\em loc. cit.} for the discussion of numerous
applications of the dynamical Weyl group.
Let $U_\hbar$ be the ``graded enveloping" algebra of $\check\fg$, i.e. the graded $\BC[\hbar]$-algebra generated
by $\check\fg$ with relations $xy-yx=\hbar[x,y]$ for $x,y\in\check\fg$. It is a subalgebra of $U(\check\fg)\otimes\BC[\hbar]$. Any $V\in\on{Rep}(\check G)$ gives rise to a structure of $U_\hbar$-module
on $V\otimes\BC[\hbar]$, and then the construction of~\cite{EV} (Section~4, Main Definition~2) produces the {\em dynamical Weyl group} operators $A_{w,V}:\ V\otimes\BC(\check\ft{}^*\times\BA^1)\to
V\otimes\BC(\check\ft{}^*\times\BA^1)$. Here $\BA^1$ is the affine line with coordinate $\hbar$.
Note that for $\mu\in\Lambda$, the operator $A_{w,V}$ takes the weight space
$V_\mu\otimes\BC(\check\ft{}^*\times\BA^1)$ to the weight space $V_{w\mu}\otimes\BC(\check\ft{}^*\times\BA^1)$.

Recall that $\check\ft{}^*=\ft$, and let $2\rho\in\check\ft{}^*=\ft$ be the sum of all positive roots of $\check G$. Let $e$ stand for the neutral element of the Weyl group $W$. The following is the main result of this paper:

\th{main} Let $V\in\on{Rep}(\check G),\ \mu\in\Lambda^+$. For $x\in\ft\times\BA^1$ we have
$\sA_{e,w}^{S_\mu(V)}(x)=A_{w,V}(-x-\hbar\rho)$. This is an equality of rational functions on $\ft\times\BA^1$
with values in $\Hom(V_\mu,V_{w\mu})$.
\eth
\ssec{}{Organization of the paper}
\refs{proof} is devoted to the proof of \reft{main}. Let us note that the only proof we have at the moment is by means
of some ``brute force" calculation. In particular, we don't use the elegant definition of $A_{w,V}$ via intertwining
operators between certain Verma-type modules over $\check\fg$ given in \cite{EV};
rather by using some formal properties of both $A_{w,V}$ and $\sA_{w,y}^{S_{\mu}(V)}$ we reduce \reft{main} to the case
of $G=PGL(2)$ where the two sides are compared by means of an explicit calculation (for $A_{w,V}$ this calculation is
done in \cite{EV}). It would be interesting to find a more conceptual proof of \reft{main}.

In \refs{aff} we formulate a conjectural generalization of \reft{main} to the case of affine Kac-Moody groups
(based on the construction of affine analogs of $\ol{\CW}{}^\lambda_\mu$ given in \cite{BF1}).

\ssec{ack}{Acknowledgments}
The idea of this work was born after a series of conversations we had with D.~Maulik and A.~Okounkov, who explained
to us the contents of their work in progress \cite{MO}; in particular, one of the (numerous) results of \cite{MO} probably
implies a proof of our \refco{affi} in the case of level 1 representations for $G=GL(n)$.
In fact, one of the main purposes of this paper was to put
this result in some more general framework. We would like to thank D.~Maulik and A.~Okounkov for their
patient explanations of the contents of \cite{MO} to us.

We are also grateful to P.~Etingof for teaching us about the dynamical Weyl group.
We are happy to thank the IAS at the Hebrew University of Jerusalem for the excellent working conditions.
M.~F. was partially supported by the RFBR grant 09-01-00242, the Ministry of Education and
Science of Russian Federation, grant No. 2010-1.3.1-111-017-029,
and the AG Laboratory HSE, RF government grant, ag. 11.G34.31.0023.
\sec{proof}{Proofs}

\ssec{prlem}{Proof of~\refl{ff}}

(a) By definition of the Mirkovi\'c-Vilonen fiber functor, we have $r^e_\mu S(V)=V_\mu$. Due to the $G(\CO)$-equivariance of the sheaf $S(V)$, we have $r^w_\mu S(V)\cong r^e_{w\mu}S(V)=V_{w\mu}$.
Let $\lambda\in\Lambda^+$ be big enough so that $\ol\Gr{}^\lambda$ contains the support of $S(V)$.
According to~Propositions~1.3.1 and~1.3.2 of~\cite{KT}, there is a $T\times\BC^*$-invariant open subset
$U^\mu\subset\ol\Gr{}^\mu$ containing $\mu$ and an open $T\times\BC^*$-equivariant embedding
$U^\lambda_\mu:=U^\mu\times\ol\CW{}^\lambda_\mu\hookrightarrow\ol\Gr{}^\lambda$.

Recall that $2\rho$ stands for the sum of all positive roots of $\check G$. We view $2\rho$ as a
cocharacter $\BC^*\to T$. Then $\fT^w_\mu$ is the repellent of the point $\mu\in\Gr$ with respect to the one-parametric subgroup $2w\rho:\ \BC^*\to T$, i.e. $\fT^w_\mu=\{g\in\Gr:\ \lim_{z\to\infty}2w\rho(z)g=\mu\}$
(see~\cite{MV},~(3.5)).
It follows that $\fT^w_\mu\cap U^\lambda_\mu=(\fT^w_\mu\cap U^\mu)\times(\fT^w_\mu\cap\ol\CW{}^\lambda_\mu)$.
Now $S(V)|_{U^\lambda_\mu}\cong\IC(\ol\Gr{}^\mu)|_{U^\mu}\boxtimes S_\mu(V)$. We conclude that
$r^w_\mu S(V)\cong r^w_\mu\IC(\ol\Gr{}^\mu)\otimes\sr^w_\mu S_\mu(V)$.

Now $\fT^w_\mu\cap\ol\Gr{}^\mu=N_-^w(\CO)\cdot\mu$ (see~\cite{MV}, proof of~Theorem~3.2) is an affine space
contained in $\Gr^\mu$. So $\IC(\ol\Gr{}^\mu)$ is constant along $\fT^w_\mu\cap\ol\Gr{}^\mu$, and we have
canonically $r^w_\mu\IC(\ol\Gr{}^\mu)\cong\BC$. Thus, $V_{w\mu}\cong r^w_\mu S(V)\cong\BC\otimes\sr^w_\mu S_\mu(V)=\sr^w_\mu S_\mu(V)$. The proof of (a) is complete.

To prove (b) note that by (a) and~Theorem~3.5 of~\cite{MV}, $\sr^w_\mu S_\mu(V)$ is concentrated in degree 0 and
coincides with the degree zero part of $\sR^w_\mu S_\mu(V)$. Thus the natural map (the action)
$\sr^w_\mu S_\mu(V)\otimes H^\bullet_{T\times\BC^*}(pt)\to\sR^w_\mu S_\mu(V)$ is the desired isomorphism.
The proof of the lemma is complete.

\ssec{prtr}{Proof of~\reft{main} for $G=PGL(2)$}

For $G=PGL(2)$ we have $\check G=SL(2)$, and let $V=V^\lambda$ be the irreducible representation of $\check G$ with
highest weight $\lambda\in\BN$, and highest vector $v_\lambda$. We consider the basis $\{v_\mu,\ \mu=-\lambda, -\lambda+2,\ldots,\lambda-2,\lambda\}$ of $V^\lambda$ such that $v_{\lambda-2k}=\frac{f^k}{k!}v_\lambda$ where
$(e,h,f)$ is the standard basis of $\fsl(2)$. In this basis, the action of dynamical Weyl group is computed in
Propositions~6 and~12 of~\cite{EV}. If $s$ is the nontrivial element of $W=\BZ/2\BZ$, and $\mu\geq0$, then
\eq{eti}
A_{s,V^\lambda}(x)v_\mu=(-1)^{\frac{\lambda-\mu}{2}}
\frac{(x+2\hbar)(x+3\hbar)\ldots(x+\hbar+\frac{\lambda-\mu}{2}\hbar)}
{(x+\hbar-\frac{\lambda+\mu}{2}\hbar)(x+2\hbar-\frac{\lambda+\mu}{2}\hbar)\ldots(x-\mu\hbar)}v_{-\mu}
\end{equation}
Here $x\in\BC[\ft]$ stands for the linear function $\langle\cdot,\check\alpha\rangle$ where $\check\alpha$ is the positive simple root of $G$. We see that
\eq{etin}
A_{s,V^\lambda}(-x-\hbar)v_\mu=
\frac{(x-\hbar)(x-2\hbar)\ldots(x-\frac{\lambda-\mu}{2}\hbar)}
{(-x-\frac{\lambda+\mu}{2}\hbar)(-x+\hbar-\frac{\lambda+\mu}{2}\hbar)\ldots(-x-\hbar-\mu\hbar)}v_{-\mu}
\end{equation}

We have $S(V^\lambda)=\unl\BC{}_{\ol\Gr{}^\lambda}[\lambda]$, and $V^\lambda=H^\bullet(\ol\Gr{}^\lambda,\unl\BC[\lambda])$. The basis $\{v_\mu\}$ consists of the (Poincar\'e duals) of the fundamental classes of the closures of intersections $\fT_\mu^e\cap\ol\Gr{}^\lambda$ (see e.g.~Section~5.2 of~\cite{BF}). We conclude from the proof of~\refl{ff} that under the isomorphism
$\sR^e_\mu S_\mu(V^\lambda)\cong V^\lambda_\mu\otimes\BC[x,\hbar]$
(resp. $\sR^s_\mu S_\mu(V^\lambda)\cong V^\lambda_{-\mu}\otimes\BC[x,\hbar]$) the positive generator of the
cohomology with integral coefficients goes to $v_\mu$ (resp. to $v_{-\mu}$). Here $x$ stands for the positive
generator of the integral equivariant cohomology $H^2_T(pt,\BZ)$.
Let $1_\mu$ stand for the unit element in the equivariant costalk $\mu^!\unl\BC{}_{\ol\Gr{}^\lambda}[\lambda]$.

Then
\eq{odin}
C^e_\mu(1_\mu)=(-x-\frac{\lambda+\mu}{2}\hbar)(-x+\hbar-\frac{\lambda+\mu}{2}\hbar)\ldots(-x-\hbar-\mu\hbar)v_\mu,
\end{equation}
while
\eq{dva}
C^s_\mu(1_\mu)=(x-\hbar)(x-2\hbar)\ldots(x-\frac{\lambda-\mu}{2}\hbar)v_{-\mu}.
\end{equation}
In effect, $\fT^e_\mu\cap\ol\CW{}^\lambda_\mu$ is a $T\times\BC^*$-invariant open subset of the affine space
$\BA^{\frac{\lambda-\mu}{2}}\simeq\fT^e_\mu\cap\ol\Gr{}^\lambda$ containing the only $T\times\BC^*$-fixed
point $\mu$, and the action of $T\times\BC^*$ on $\BA^{\frac{\lambda-\mu}{2}}$ is linear with weights
$\{-x-\frac{\lambda+\mu}{2}\hbar,\ -x+\hbar-\frac{\lambda+\mu}{2}\hbar,\ldots,\ -x-\hbar-\mu\hbar\}$
(see e.g.~Section~5.2 of~\cite{BF}). Similarly, $\fT^s_\mu\cap\ol\CW{}^\lambda_\mu$ is a $T\times\BC^*$-invariant open subset of the affine space $\BA^{\frac{\lambda-\mu}{2}}$ with the linear $T\times\BC^*$-action with weights
$\{x-\hbar,\ x-2\hbar,\ldots,\ x-\frac{\lambda-\mu}{2}\hbar\}$.

Comparing~\refe{odin} and~\refe{dva} with~\refe{etin} we conclude that $\sA^{S_\mu(V^\lambda)}_{e,s}(x)=
A_{s,V^\lambda}(-x-\hbar)$.

\ssec{red}{Reduction to rank 1}

By construction of operators $\sA$, for any $\mu\in\Lambda,\ V\in\on{Rep}(\check G),\ w_1,w_2\in W$ we have
$\sA^{S_\mu(V)}_{w_2,w_1w_2}\circ\sA^{S_\mu(V)}_{e,w_2}=\sA^{S_\mu(V)}_{e,w_1w_2}$, and also
$\sA^{S_\mu(V)}_{w_2,w_1w_2}(x)=\sA^{S_\mu(V)}_{e,w_1}(w_2x)$. In particular,
\eq{tri}
\sA^{S_\mu(V)}_{e,w_1}(w_2x)\circ\sA^{S_\mu(V)}_{e,w_2}(x)=\sA^{S_\mu(V)}_{e,w_1w_2}(x).
\end{equation}
According to~Lemma~17 of~\cite{EV}, in case length $l(w_1w_2)=l(w_1)+l(w_2)$, we have
\eq{chetyre}
A_{w_1,V}(w_2x)\circ A_{w_2,V}(x)=A_{w_1w_2,V}(x).
\end{equation}
Comparing~\refe{tri} and~\refe{chetyre} and choosing a reduced decomposition of $w\in W$ into a product of
simple reflections, we see that it is enough to prove~\reft{main} for $w=s$ a simple reflection.

Let $G\supset P\twoheadrightarrow L$ be the corresponding (to $s$) submimimal parabolic subgroup containing $B_-$, and Levi group. We denote by $\check L\subset\check G$ the corresponding semisimple rank 1 Levi subgroup of $\check G$. The connected components of $\Gr_P$ are numbered by the lattice $\Lambda_{G,P}$ of characters of
the center $Z(\check L)$ of $\check L$. For $\theta\in\Lambda_{G,P}$ we denote by $\imath_\theta:\ \Gr_{P,\theta}\hookrightarrow\Gr_G$ the locally closed embedding of the corresponding connected component.
We denote by $\fp^\theta:\ \Gr_{P,\theta}\to\Gr_{L,\theta}$ the projection to the corresponding connected
component of the affine Grassmannian of $L$. Let $\check\rho_{G,L}$ stand for $\check\rho_G-\check\rho_L$.
Then for a $G(\CO)$-equivariant perverse sheaf $\CF$ on $\Gr_G$ the direct sum $\bigoplus_{\theta\in\Lambda_{G,P}}\fp^\theta_*\imath_\theta^!\CF[\langle\mu,2\check\rho_{G,L}\rangle]$ is a
$L(\CO)$-equivariant perverse sheaf on $\Gr_L$ (see~\cite{BD},~Proposition~5.3.29) to be denoted by $R^{G,L}\CF=\bigoplus_{\theta\in\Lambda_{G,P}}R^{G,L}_\theta\CF$. Moreover, according to {\em loc. cit.},
for any $V\in\on{Rep}(\check G)$ we have $R^{G,L}S(V)\cong S_L(\on{Res}^{\check G}_{\check L}V)$ where
$\on{Res}^{\check G}_{\check L}V$ stands for the restriction of $V$ from $\check G$ to $\check L$.

Given $\mu\in\Lambda$ we will denote by $\bar\mu$ its image in the quotient lattice $\Lambda_{G,P}$.
For a $T\times\BC^*$-equivariant constructible complex on the slice $\CW_\mu$, it is easy to see that
$\sR^{G,L}_{\bar\mu}\CF=\fp^{\bar\mu}_*\imath_{\bar\mu}^!\CF$ will be a $T\times\BC^*$-equivariant constructible
complex supported on $\CW_\mu^L$ (the transversal slice on the affine Grassmannian of $L$ through the point $\mu$). For $V\in\on{Rep}(\check G),\ \mu\in\Lambda^+$, and $w\in W_L=\{e,s\}$ by the base change we have $\sR^{G,w}_\mu S^G_\mu(V)=\sR^w_\mu S_\mu(V)\cong \sR^{L,w}_\mu \sR^{G,L}_{\bar\mu}S^G_\mu(V)\cong\sR^{L,w}_\mu S^L_\mu(\on{Res}^{\check G}_{\check L}V)$ where we have added the superscripts $^G$ or $^L$ to $\sR^w_\mu$ and $S_\mu$ in order
to distinguish the hyperbolic stalks taken on transversal slices of the affine Grassmannians of $G$ of $L$.
Moreover, $\sA_{e,s}^{G,S^G_\mu(V)}=\sA_{e,s}^{L,S^L_\mu(\on{Res}^{\check G}_{\check L}V)}$.

Now all the connected components of $\Gr_L$ are of the type considered in~\refss{prtr}, and the $L(\CO)$-equivariant perverse sheaf $S^L_\mu(\on{Res}^{\check G}_{\check L}V)$ is a direct sum of sheaves considered in~\refss{prtr}. Finally, the action of the dynamical Weyl group operator $A_{s,V}$ is defined
as $A_{s,\on{Res}^{\check G}_{\check L}V}$. By the virtue of~\refss{prtr} we conclude that
$\sA^{S_\mu(V)}_{e,s}(x)=A_{s,V}(-x-\hbar\rho)$.

This completes the proof of~\reft{main}.

\sec{aff}{Affine case}

\ssec{gaff}{The group $G_{\aff}$}From now on we assume that $G$ is almost simple and simply connected.
To a connected reductive group $G$ as above one can associate the corresponding
affine Kac-Moody group $G_{\aff}$ in the following way.
One can consider the polynomial loop group $G[t,t^{-1}]$ (this is an infinite-dimensional group ind-scheme)

It is well-known that $G[t,t^{-1}]$ possesses a canonical central extension $\tilG$ of $G[t,t^{-1}]$:
$$
1\to \GG_m\to \tilG\to G[t,t^{-1}]\to 1.
$$
Moreover,\ $\tilG$ has again a natural structure of a group ind-scheme.

The multiplicative group $\GG_m$ acts naturally on $G[t,t^{-1}]$ and this action lifts to $\tilG$.
We denote the corresponding semi-direct product by $G_{\aff}$; we also let $\grg_{\aff}$ denote its Lie algebra.

The Lie algebra $\grg_{\aff}$ is an untwisted affine Kac-Moody Lie algebra.
In particular,
it can be described by the corresponding affine root system. We denote by $\grg_{\aff}^{\vee}$ the
{\em Langlands dual affine Lie algebra} (which corresponds to the dual affine root system)
and by $G^{\vee}_{\aff}$ the corresponding dual affine Kac-Moody group, normalized by the property
that it contains $G^{\vee}$ as a subgroup (cf. \cite{BF1}, Subsection 3.1 for more details).

We denote by $\Lam_{\aff}=\ZZ\x\Lam\x\ZZ$ the coweight lattice of $G_{\aff}$; this is the same as the
weight lattice of $G_{\aff}^{\vee}$. Here the first $\ZZ$-factor
is responsible for the center of $G_{\aff}^{\vee}$ (or $\hatG^{\vee}$);
it can also be thought of as coming from the loop
rotation in $G_{\aff}$. The second $\ZZ$-factor is responsible for the loop rotation in $G_{\aff}^{\vee}$  it may also be thought of
as coming from the center of $G_{\aff}$).
We denote by $\Lam_{\aff}^+$ the set of dominant weights of $G_{\aff}^{\vee}$ (which is the same as the set of dominant
coweights of $G_{\aff}$). We also denote by $\Lam_{\aff,k}$ the set of weights of $G_{\aff}^{\vee}$ of level $k$,
i.e. all the weights of the form $(k,\olam,n)$. We put $\Lam_{\aff,k}^+=\Lam_{\aff}^+\cap \Lam_{\aff,k}$.

\smallskip
\noindent
{\bf Important notational convention:}
From now on we shall denote elements
of $\Lam$ by $\blambda,\bmu...$ (instead of just writing $\lam,\mu...$ in order to distinguish them from
the coweights of $G_{\aff}$ (= weights of $G_{\aff}^{\vee}$), which we shall just denote by
$\lam,\mu...$

Let
$\Lam_k^+\subset \Lam$ denote the set of dominant coweights of $G$ such that $\la \blambda,\alp)\leq k$
when $\alp$ is the highest root of $\grg$.
Then it is well-known that a weight $(k,\olam,n)$ of
$G_{\aff}^{\vee}$ lies in $\Lam_{\aff,k}^+$  if and only if $\olam\in\Lam_k^+$ (thus $\Lam_{\aff,k}=\Lam_k^+\x \ZZ$).

Let also $W_{\aff}$ denote affine Weyl group of $G$ which is the semi-direct product of $W$ and $\Lam$.
It acts on the lattice $\Lam_{\aff}$ (resp. $\hatLam$)
preserving each $\Lam_{\aff,k}$ (resp. each $\hatLam_k$). In order to describe this action explicitly it is convenient
to set $W_{\aff,k}=W\ltimes k\Lam$ which naturally acts on $\Lam$. Of course the groups $W_{\aff,k}$ are canonically
isomorphic to $W_{\aff}$ for all $k$. Then the restriction of the $W_{\aff}$-action to $\Lam_{\aff,k}\simeq\Lam\x\ZZ$
comes from the natural $W_{\aff,k}$-action on the first multiple.

It is well known that every $W_{\aff}$-orbit on $\Lam_{\aff,k}$ contains unique element of $\Lam_{\aff,k}^+$.
This is equivalent to saying that $\Lam_k^+\simeq \Lam/W_{\aff,k}$.

\ssec{rev3}{Repellents}
Our dream is to create an analog of the affine Grassmannian $\Gr_G$ and the above results
about it in the case when $G$ is replaced by the (infinite-dimensional) group $G_{\aff}$.
The first attempt to do so was made in~\cite{BF1}: namely, in {\em loc. cit.} we have constructed analogs $\ol\CW{}^\lambda_\mu$ of the varieties $\ol\CW{}^{\blambda}_{\bmu}$ in the case when $G$ is replaced by  $G_{\aff}$.
We will not reproduce the (rather involved) definition of $\ol\CW{}^\lambda_\mu$ here; we only mention that the
open pieces $\CW{}^\lambda_\mu$ are the moduli spaces of certain $G$-bundles on $\BP^2$ trivialized at infinity $\BP^1\subset\BP^2$. The complement $\BA^2=\BP^2\setminus\BP^1$ is equipped with coordinates $(z,t)$, and with
the action of $\BC^*\times\BC^*:\ (a,b)\cdot(z,t)=(az,bt)$. The 2-dimensional torus $\BC^*\times\BC^*$ acts on
$\ol\CW{}^\lambda_\mu$ by transport of structure, and $T\subset G$ acts via trivialization at infinity.
Let $\BC^*_{\on{hyp}}:=\{(c,c^{-1})\}\subset\BC^*\times\BC^*$ (resp.  $\BC^*_\Delta:=\{(c,c)\}\subset\BC^*\times\BC^*$) stand for the antidiagonal, alias hyperbolic, (resp. diagonal) subgroup. Let us denote the torus $\BC^*_{\on{hyp}}\times T$ by $\widehat T$. Thus, the slices $\ol\CW{}^\lambda_\mu$ are equipped with the action of the torus $\widehat{T}\times\BC^*_\Delta$.
The Lie algebra of $\widehat{T}$ is denoted by $\widehat\ft$.

The cocharacter lattice $X_*(\widehat{T})\cong\BZ\times\Lambda\subset\BZ\times\Lambda\times\BZ=X_*(T_{\aff})$.
Accordingly, we have $\widehat{T}\subset T_{\aff}$. We propose that the action of $\widehat{T}$ on $\ol\CW{}^\lambda_\mu$ of the previous paragraph is nothing else than the action coming from the $T_{\aff}$-action on the transversal slices of the double affine Grassmannian. Also, the action of $\BC^*_\Delta$ on $\ol\CW{}^\lambda_\mu$ of the previous paragraph is nothing but the loop rotation action on the transversal slices of the double affine Grassmannian.

Let $I$ (resp. $I_{\aff}=I\sqcup i_0$) stand for the set of vertices of the Dynkin diagram of $G$ (resp. of $G_{\aff}$). For $i\in I$ we denote by $\ol\omega_i\in X_*(T)$ the corresponding fundamental coweight of $G$,
and we denote by $a_i\in\BN$ the corresponding label of the Dynkin diagram of $G_{\aff}$.
Then $\omega_{i_0}:=(1,0)\in\BZ\times X_*(T)=X_*(\widehat{T}),\ \omega_i:=(a_i,\ol\omega_i)\in X_*(\widehat{T}),\ i\in I$, are the fundamental coweights of $G_{\aff}$. We set $\rho:=\sum_{i\in I_{\aff}}\omega_i$ (not to be confused with the halfsum $\bar\rho$ of positive coroots of $G$). The group $W_{\aff,k}$ acts on $X_*(\widehat{T})$, and for $w\in W_{\aff,k}$ let us view $w\rho$ as a one-parametric subgroup $\BC^*\to\widehat T$. The torus $\widehat T$ acts on $\ol\CW{}^\lambda_\mu$ with the only fixed point, to be denoted abusively by $\mu$, and by analogy with~\refss{prlem}, we define $\fT^w_\mu\subset\ol\CW{}^\lambda_\mu$ as the repellent
$\fT^w_\mu:=\{g\in\ol\CW{}^\lambda_\mu:\ \lim_{c\to\infty}w\rho(c)g=\mu\}$. For $G=SL(N)$ these repellents were
considered by H.~Nakajima in~Section~6 of~\cite{Na} under the name of MV cycles.

The same procedure as in~\refss{hyp} defines an isomorphism $\sA^{\IC(\ol\CW{}^\lambda_\mu)}_{w,y}:\
\sR^w_\mu\IC(\ol\CW{}^\lambda_\mu)\otimes_{\BC[\widehat\ft\times\BA^1]}\BC(\widehat\ft\times\BA^1)\iso
\sR^y_\mu\IC(\ol\CW{}^\lambda_\mu)\otimes_{\BC[\widehat\ft\times\BA^1]}\BC(\widehat\ft\times\BA^1)$. Here $\BA^1$ is the affine line with coordinate $\hbar$, the positive generator of the integral equivariant cohomology
$H^2_{\BC^*_\Delta}(pt,\BZ)$.

\conj{affi}
(a) For $\lambda,\mu\in\Lambda^+_k,\ w\in W_{\aff,k}$, we have $\sR^w_\mu\IC(\ol\CW{}^\lambda_\mu)\cong V^\lambda_{w\mu}\otimes\BC[\widehat\ft\times\BA^1]$ where
$V^\lambda_\mu$ stands for the $w\mu$-weight space of the irreducible integrable $G^\vee_{\aff}$-module with highest weight $\lambda$;

(b) We have $\sA^{\IC(\ol\CW{}^\lambda_\mu)}_{e,w}(x)=A_{w,V^\lambda}(-x-\hbar\rho)$ where $A_{w,V^\lambda}$ stands for the dynamical affine Weyl group action of~\cite{EV}.
\econj

As was mentioned in the Introduction, in the case of representations of level 1, one can probably deduce \refco{affi} for $G=SL(n)$ (as well as some generalization
of \refco{affi} to $G=GL(n)$) from \cite{MO} (by using also the results of \cite{BF1} and \cite{Na}).

\bigskip
\footnotesize{
{\bf A.B.}: Department of Mathematics, Brown University,
151 Thayer St., Providence RI
02912, USA;\\
{\tt braval@math.brown.edu}}

\footnotesize{
{\bf M.F.}: IMU, IITP and State University Higher School of Economics\\
Department of Mathematics, 20 Myasnitskaya st, Moscow 101000 Russia;\\
{\tt fnklberg@gmail.com}}


\begin{thebibliography}{jjjj}
\bibitem{AST} R.~Asherova, Yu.~Smirnov, V.~Tolstoy, {\em Projection operators for simple Lie groups}, Theor. Math. Phys. {\bf 8} (1971), 813--825.

\bibitem{BD}
A.~Beilinson and V.~Drinfeld, {\em Quantization of Hitchin's Hamiltonians
and Hecke eigen--sheaves}, Preprint, available at
http://www.math.uchicago.edu/$\widetilde{\hphantom{m}}$mitya/langlands.html

\bibitem{BF} R.~Bezrukavnikov, M.~Finkelberg, {\em Equivariant Satake Category and Kostant-Whittaker Reduction},
Moscow Math. Jour. {\bf 8} (2008), 39--72.


\bibitem{BF1}
A.~Braverman and M.~Finkelberg, {\em Pursuing the double affine Grassmannian I:
Transversal slices via instantons on $A_k$-singularities},
Duke Math. J. {\bf 152} (2010), 175--206.






\bibitem{EV} P.~Etingof, A.~Varchenko, {\em Dynamical Weyl groups and Applications},
Adv. Math. {\bf 167} (2002), 74--127.

\bibitem{Gi}
V.~Ginzburg, {\em Perverse sheaves on a loop group and Langlands
duality}, Preprint, alg--geom/9511007.

\bibitem{KT}
M.~Kashiwara, T.~Tanisaki, {\em Kazhdan-Lusztig conjecture for affine Lie algebras with negative level},
Duke Math. J. {\bf 77} (1995), 21--62.

\bibitem{Lu-qan}
G.~Lusztig,
{\em Singularities, character formulas, and a $q$-analog of
weight multiplicities}, Analysis and topology on singular spaces,
Ast{\'e}risque, {\bf 101-102} (1983), 208--229.

\bibitem{MO}
D.~Maulik and A.~Okounkov, in preparation.

\bibitem{MV}
I.~Mirkovi\'c and K.~Vilonen, {\em Geometric Langlands duality and
  representations of algebraic groups over commutative rings},
Annals of Math. (2) {\bf 166} (2007), 95--143.

\bibitem{Na}
H.~Nakajima, {\em Quiver varieties and branching}, SIGMA {\bf 5}
(2009), 003, 37 pages.

\bibitem{Z} D.~P.~Zhelobenko, {\em Extremal projectors and generalized Mickelsson algebras on reductive Lie algebras}, Math. USSR, Izv. {\bf 33} (1989), 85--100.


\end{thebibliography}
\end{document}